\documentclass[12pt]{amsart}
\usepackage[left=1in, right=1in, top=1in, bottom =1in]{geometry}

\usepackage{settings} % add packages, settings, and declarations in settings.tex

\begin{document}
%%%%%%%%%%%%%%%%%%%%%%%%%%%%%%%%%%%%%%
%%%%%%%%% TITLE + AUTHORS %%%%%%%%%%%%
%%%%%%%%%%%%%%%%%%%%%%%%%%%%%%%%%%%%%%
\title{Optimal land conservation decisions for multiple species}

\author[Buhler]{Cassidy K. Buhler}
\address{Department of Decision Sciences and MIS, Bennett S. LeBow College of Business, Drexel University, Philadelphia, PA 19104}
\email{cb3452@drexel.edu}

\author[Benson]{Hande Y. Benson}
\address{Department of Decision Sciences and MIS, Bennett S. LeBow College of Business, Drexel University, Philadelphia, PA 19104}
\email{hvb22@drexel.edu}

%%%%%%%%%%%%%%%%%%%%%%%%%%%%%%%%%%%%%%
%%%%%%%%%%%%% ABSTRACT %%%%%%%%%%%%%%%
%%%%%%%%%%%%%%%%%%%%%%%%%%%%%%%%%%%%%%
\begin{abstract}
Given an allotment of land divided into parcels, government decision-makers, private developers, and conservation biologists can collaborate to select which parcels to protect, in order to accomplish sustainable ecological goals with various constraints. In this paper, we propose a mixed-integer optimization model that considers the presence of multiple species on these parcels, subject to predator-prey relationships and crowding effects.
\end{abstract}

\maketitle %Puts title on a title page (except if this is an article, then it just goes at the top)

\section{Introduction} %%%%%%%%%%%%%%%
% Ecological models are used to estimate spatial population of animals in their habitats. 

% Preserving habitats for wildlife species is imperative to biodiversity, and ultimately, human survival. Ironically, the loss of these habitats are largely driven by human activities.
\subsection{Motivation}

Climate change poses a  massive threat to the health of humanity \cite{pruss2016preventing}. The ramifications of the climate crisis extend beyond increased global temperatures, as necessities like food and water will become more and more scarce. Humans are not the only species that will be affected. At this rate, it's been estimated that one in six species could face extinction \cite{urban2015accelerating}. In an interconnected world, it is, therefore, important to incorporate sustainability with a focus on biodiversity into every level of decision making.

 Government decision-makers at the municipal, county, state, and federal levels frequently work with private companies, including engineering firms and developers, to determine where development should take place and where ``protected areas'' must be established.  That is, spaces where species are protected from human interference. In this paper, we will be focusing on the design of protected areas.

\subsection{Designing protected areas}
Given an allotment of land divided into parcels, our task is to select which parcels to protect, in order to accomplish ecological goals subject to various financial and economic constraints.  These decisions are aided by optimization models. \cite{alagador2022operations, billionnet2013mathematical}

Protected areas can take on many configurations. Recent work in this area is developing models to ensure particular spatial properties. For instance, some researchers have focused on selecting parcels such that the protected area is connected \cite{costanza2020preserving, gupta2019reserve, dilkina2017trade} or contiguous \cite{wang2020optimizing, onal2016optimal}. However, many models that include these additional constraints focus on a single species; incorporating multiple species adds more dimensions to an already computationally expensive problem to solve.

% In this paper, we present a simple model for multiple species.

% In this paper, we address an existing model \cite{gupta2019reserve} and extend it towards multiple species. 

% \subsection{Gupta's paper - not sure where to add this in}
% Gupta's paper uses a spatial capture-recapture model (SCR) to estimate spatial population of a single species in particular area. Using this, they maximize the protected density weighted-connectivity \cite{morin2017model, sutherland2015modelling} of a reserve and applied home-range constraints.

\section{Model and Methodologies} %%%%%%%%%%%%%%

In this section, we present the optimization models for making land preservation decisions in the presence of multiple species and budget constraints.  We start by introducing our notation, then present a model where the species do not interact and one where predator-prey relationships are present.

\subsection{Notation}
The notation used in the optimization models and in the rest of the paper are as follows:

\vspace{0.2in}

Sets:

\begin{tabular}{ll}
$P$		& 	the set of parcels \\
$S$		& 	the set of species
\end{tabular}

\vspace{0.2in}

Parameters:

\begin{tabular}{lp{4.5in}}
$w_i$		&	Weight to prioritize species $i$ \\
$N_i(p)$	&	Number of individuals of species $i$ observed at parcel $p$ \\
$N_i$       &   Total population of species $i$ across all parcels \\
$\tilde{N}_i(p)$ &	Number of individuals of species $i$ that are simulated to be at parcel $p$ in the future \\
$c_p$		&	The cost of preserving parcel $p$ \\
$B$		&	Budget 
\end{tabular}

\vspace{0.2in}

Decision Variables:

\begin{tabular}{ll}
$x_p$		& 	Binary variable denoting whether or not a parcel is preserved
\end{tabular}

\subsection{Model without Species Interaction}

Our baseline model assumes that land preservation decisions can be made by taking into account only the number of individuals of each species as currently observed.  This means we assume the different species do not interact with each other and that crowding effects do not occur.  Another interpretation is that the model takes only the present conditions into consideration, instead of focusing on sustainability.

This baseline model \eqref{withoutInteractionModel} is a knapsack problem, where the objective is to save a weighted combination of species' populations subject to a budget constraint.  The weights, $w_i$, can be chosen to prioritize endangered species or reflect other conservation concerns.  
\begin{equation} 
	\begin{array}{ll}
	\text{maximize } & \displaystyle \sum_{i \in S} w_i \sum_{p \in P} N_i(p) x_p \\
	\text{subject to } & \displaystyle \sum_{p \in P} c_p x_p \le B \\
	& x_p \in \{0,1\} \quad \forall p \in P
         \label{withoutInteractionModel}
	\end{array}
\end{equation}

\subsection{Model with Species Interaction}
The main model we present is the model with species interaction, such as predator-prey relationships.  The model itself has the same overall form as \eqref{withoutInteractionModel} with one critical difference.  The species' populations $\tilde{N}_i(p)$ are calculated using a simulation that models the numbers of individuals present in each parcel after $T$ time periods.  In our numerical testing, we have taken $T = 2000$, which represents steady-state populations.
\begin{equation} 
	\begin{array}{ll}
	\text{maximize } & \displaystyle \sum_i w_i \sum_{p} \tilde{N}_i(p) x_p \\
	\text{subject to } & \displaystyle \sum_p c_p x_p \le B \\
	& x_p \in \{0,1\} \quad \forall p \in P
         \label{withInteractionModel}
	\end{array}
\end{equation}

$\tilde{N}_i(p)$ is obtained from $N_i(p)$ using the Gavina et.al.'s model from \cite{gavina2018multi}, which 
adapts the classical Lotka-Volterra equations describing predator-prey relationships to multiple species and 
takes into account crowding effects.  The simulation is described in Algorithm 1.

\begin{algorithm}
\caption{Lotka-Volterra competition with crowding effects}\label{alg:LV-crowding}
	\SetKwInOut{Input}{input}\SetKwInOut{Output}{output}
	\Input{$N_i(p)$: Number of individuals of species $i$ observed at parcel $p$ for each $i$ and $p$}
	\Output{$\tilde{N}_i(p)$: Number of individuals of species $i$ that are simulated to be at parcel $p$ in the future for each $i$ and $p$}
	\BlankLine
	\ForEach{parcel $p$}{
		Initialize the species counts in parcel $p$, i.e. let $\tilde{N}_i(p) = N_i(p)$ for each $i$\;
		\For{$t \leftarrow 1$ \KwTo $T$}{
			Let $\Delta N = $ (Birth Rate of species $i$)($\tilde{N}_i(p)$) - (Competition Effect of all species on $i$) - (Crowding Effect for species $i$)\;
			Let $\tilde{N}_i(p) = \tilde{N}_i(p) + \Delta N$.
		}
	}
\end{algorithm}

\section{Data} %%%%%%%%%%%%%%%%%%%%%%%

The data was generated using the code provided in Gupta et.al.'s paper \cite{gupta2019reserve}. The process is outlined below. 

A landscape is an $n\times n$ grid of parcels, where each parcel is a piece of land that can be protected. Each parcel has a value between $[0,1]$ where a higher value represents a worse habitat for that species. A landscape where better habitat parcels are clustered with each other (and the worse habitat parcels are clustered with each other) has a low habitat fragmentation. A high habitat fragmentation is an absence of these clusters, as the better habitat parcels are dispersed among worse habitat parcels, and vice versa. Based on these landscapes, density of each species are simulated. Given a species population $N_i$, these individuals are distributed among the landscape using an inhomogeneous point process. This distribution is $N_i(p)$. 

Following the model given in \cite{gavina2018multi}, we obtained $\Tilde{N}_i(p)$ by inputting $N_i(p)$ into Algorithm \ref{alg:LV-crowding}. Due to the scaling in this code, the values outputted are non-integers. To address this, we rounded the output,  $\Tilde{N}_i(p)$, to the nearest integer. 
% Then, we used the spatial capture-recapture model (SCR) \cite{sutherland2015modelling} to simulate the movement of these individuals. The SCR model is used in ecology to estimate spatial population of a single species in particular area.

% EAch species have their own habitat preferences, which we captured by assigning landscapes and populations to obtain a unique distributions of each species on a grid. 

% Each species has their own habitat preferences, which we captured by assigning them a unique distributions.  
% . Each species has their own habitat preferences, so we assigned a unique landscape to each species in order to get a distribution of that species on a grid. 

Habitat preferences were addressed by assigning landscapes and populations to species, which yields unique distributions of each species on a grid.
We generated 10,000 $10\times10$ landscapes with random habitat fragmentation levels and the two most and two least fragmented landscapes were selected. For each of these four landscapes, we explored two population sizes,  $N_i= 100$ and $ N_i = 250$. This yields 8 different distributions that represent 8 species. Additional details can be found in Table \ref{tab:speciesDeets}.

 % As outlined above, a landscape gives us a distribution and movement of individuals. Therefore, selecting multiple landscapes allows us to consider more than one species. 
% This data is inputted into our model to solve for the optimal reserve that considers the habitat preferences of each species. 

% These are the 8 different species: 
% \begin{itemize}
    % \item $S_0$ has distribution $N_0(p)$ based on highest fragmented landscape and has a population $N_0 = 100$.
%     \item  $S_1$ = high2-N100 - with distribution $N_1(p)$
%     \item  $S_2$ = high1-N250 - with distribution $N_2(p)$
%     \item  $S_3$ = high2-N250 - with distribution $N_3(p)$
%     \item  $S_4$ = low1-N100 - with distribution $N_4(p)$
%     \item  $S_5$ = low2-N100 - with distribution $N_5(p)$
%     \item  $S_6$ = low1-N250 - with distribution $N_6(p)$
%     \item  $S_7$ = low2-N250 - with distribution $N_7(p)$
% \end{itemize}

\section{Numerical Testing} %%%%%%%%%%
% We explored two landscape grid sizes: $10\times10$ and $40\times40$. 
% For each size $N$, we solved for the optimal protected area for each pair of the 4 landscapes. The scenarios are outlined in Table ~\ref{table:cases}
% Each scenario is inputted into the model for the initial time period. 

With the 8 species generated, we grouped them into sets of size $2$ and $5$ in order to implement a 2-species reserve and a 5-species reserve. For the 2-species reserves, the sets are $\{S_0, S_1\}, \{S_2, S_3\}, \{S_4, S_5\},$ and $\{S_6, S_7\}$. The 5-species reserve, the sets are $\{S_0, S_1, S_2, S_3, S_4\}$ and $\{S_5, S_6, S_7, S_0, S_1\}$. These give 6 scenarios total: 4 for 2-species reserves and 2 for 5-species reserves. 

For each scenario, we varied the budget $B$ from 0 to 100 with a stepsize of 5 and solved \eqref{withoutInteractionModel} with $N_i(p)$ and \eqref{withInteractionModel} with $\Tilde{N}_i(p)$. We tracked the similarity of each model's solution by counting the number of parcels that had the same protection status, and taking the minimum, mean, and median over all budgets except where budget is 0 and 100.  We omitted those values because the solution will always be to protect nothing or protect everything, which is uninteresting for comparison. 

\begin{figure}[ht]
    \centering
    \includegraphics[width = .8\textwidth ]{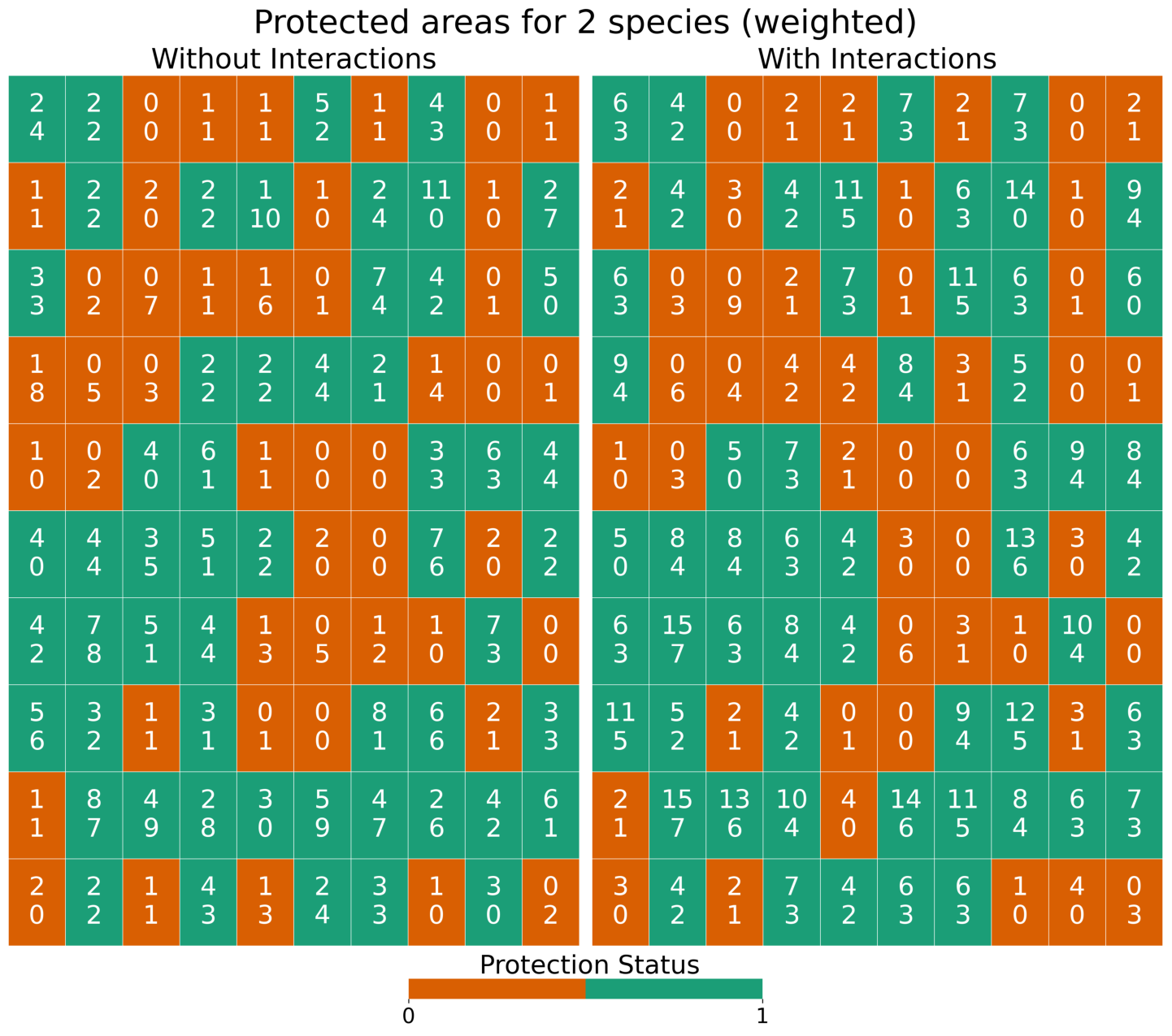}
     \caption{This graph depicts the same scenario as Figure \ref{fig:2species-b55}, but rather than setting $w_2 = w_3$, we used $w_2 = 0.9$ and $w_3 = 0.1$.  In this example, $90/100$ parcels have the same protection status.}
         \label{fig:2species-b55-weighted}
\end{figure}

\begin{figure}[ht]
    \centering
    \includegraphics[width = .95\textwidth ]{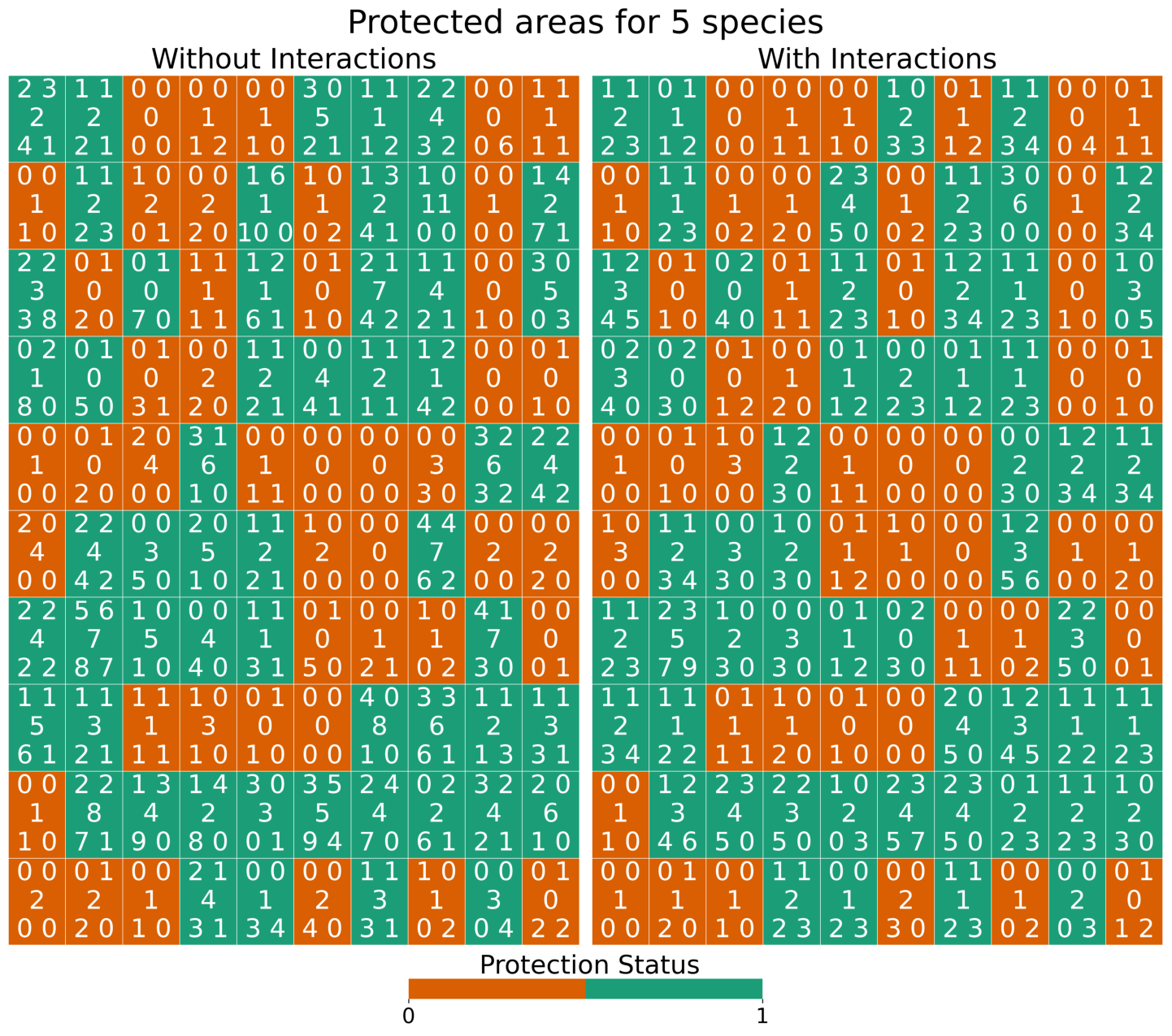}
    \caption{Displaying the solution to \eqref{withoutInteractionModel} and \eqref{withInteractionModel} for a 5-species reserve using $\{S_0, S_1, S_2, S_3, S_4\}$ when the budget is 55. Each parcel is annotated with its corresponding $N_i(p)$ or $\Tilde{N}_i(p)$ for $i = 0,\ldots,4$. The order of the values in each parcel are top left, top right, middle, bottom left, bottom right, for species $S_0, S_1, S_2, S_3, S_4$, respectively. The color of each parcel represents the preservation decision: green if preserved, orange if not. In this example, $96/100$ parcels have the same protection status.}
    \label{fig:5species-b55}
\end{figure}

\begin{figure}[ht]
    \centering
    \includegraphics[width = 0.8\textwidth ]{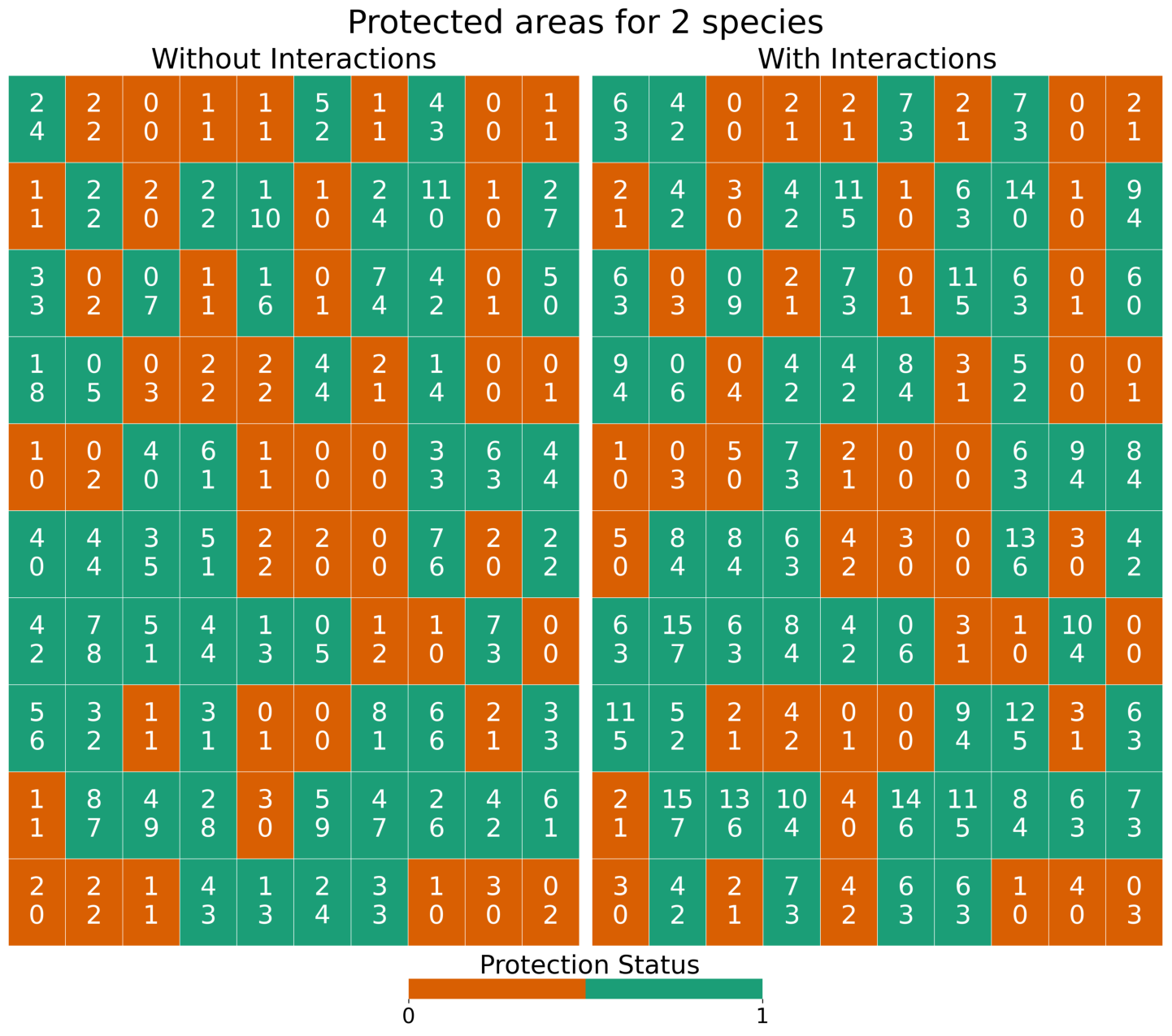}
     \caption{Displaying solutions to \eqref{withoutInteractionModel} and \eqref{withInteractionModel} for a 2-species reserve using $\{S_2, S_3\}$ and $B=55$. Each parcel is annotated with its corresponding $N_i(p)$ or $\Tilde{N}_i(p)$ for $i = 2,3$. For each parcel, the value for $S_2$ is on top and $S_3$ is on bottom. The color of each parcel represents the preservation decision: green if preserved, orange if not. In this example, $92/100$ parcels have the same protection status.}
    \label{fig:2species-b55}
\end{figure}

\begin{figure}[ht]
    \centering
    \includegraphics[width = .8\textwidth ]{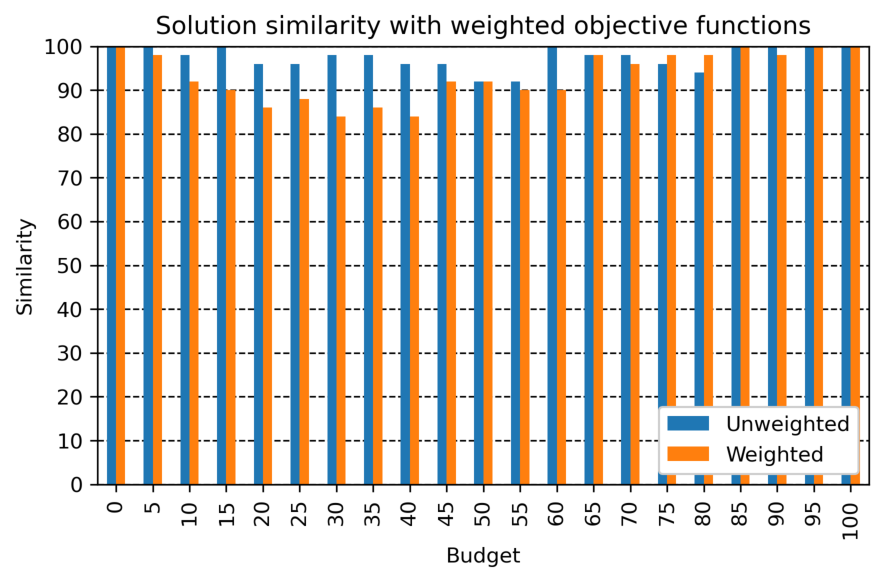}
    \caption{The similarity of solutions from \eqref{withoutInteractionModel} and \eqref{withInteractionModel} vary based on the objective function weights, $w_i$. In this example, we used $S_2$ and $S_3$ for a 2-species reserve, and solved the models where $w_2 = w_3$ (unweighted) and where $w_2 = 0.9$ and $w_3 = 0.1$ (weighted). For lower budgets, the difference is more apparent.}
    \label{fig:similarPlot}
\end{figure}

The numerical results for all cases can be found in Table \ref{tab:similar}. In addition, one solution for a 2-species reserve and 5-species reserve are included in Figure \ref{fig:2species-b55} and Figure \ref{fig:5species-b55}. The reserves found using the two models show a high degree of similarity when the objective function weights for the two species are the same. 
However, when these weights differed, the solutions showed more variation among lower budgets, which is illustrated by Figure \ref{fig:similarPlot}. A comparison between Figures \ref{fig:2species-b55} and \ref{fig:2species-b55-weighted} show a specific example of the configuration change.

\section{Discussion} %%%%%%%%%%%%%%%%%

A drawback to Algorithm \ref{alg:LV-crowding}, is that $\Tilde{N}$ needs to be rounded. It is worth exploring alternatives or variations in order to obtain integer solutions. Not only that, but it would be interesting to modify the parameters in Algorithm \ref{alg:LV-crowding} to explore dynamics that would yield a larger difference between $N$ and $\Tilde{N}$. 

For future directions, we hope to explore extensions to make the model and numerical testing more realistic. This includes increasing the grid size, expanding methods to obtain $N$, and investigating different parameters to use for the species weights and parcel costs. A larger grid size would be valuable to pursue because real-world landscapes are typically larger than $10\times10$. Also, this would allow more possible reserve solutions, thus making the problem more interesting. With regards to the methods to obtain $N$, estimating the location and movement of a species using the spatial capture-recapture (SCR) model \cite{sutherland2015modelling} and using these as inputs into our model, as done in Gupta et.al's paper, would provide a more accurate depiction of animal behavior.
% The SCR model is used in ecology to estimate spatial population of a single species in particular area. 

%%%%%%%%%%%%%%%%%%%%%%%%%%%%%%%%%%%%%%
%%%%%%%%%%%%% REFERENCES %%%%%%%%%%%%%
%%%%%%%%%%%%%%%%%%%%%%%%%%%%%%%%%%%%%%
\newpage
\bibliographystyle{plain}
\bibliography{ref.bib}

\newpage

% \newgeometry{left=.5in, right=.5in, top=.5in, bottom =.5in} 
% \setlength{\intextsep}{10pt}

\section{Appendix}
\begin{table}[ht]
    \centering
    \begin{tabular}{l|l|l}
        Species &  Fragmentation level & $N_i$\\
        \hline
      $S_0$ &  highest & 100 \\
      $S_1$ & 2nd highest & 100\\
      $S_2$ & highest & 250\\
      $S_3$ & 2nd highest & 250\\
      $S_4$ & lowest & 100\\
      $S_5$ & 2nd lowest & 100\\
      $S_6$ & lowest & 250\\
      $S_7$ & 2nd lowest & 250\\
    \end{tabular}
    \caption{Details for each species' distribution. Fragmentation level describes the density of habitat quality on a landscape, and $N_i$ is the total population of a species on that landscape. }
    \label{tab:speciesDeets}
\end{table}

\begin{table}[ht]
    \centering
    \begin{tabular}{l|l|l|l}
        Case & Min &  Average & Median \\
        \hline
         1 & 98 & 99.79  &   100\\
         2 & 92 & 97.26  &   98\\
         3 & 92 & 95.89  &   96\\
         4 & 92 & 95.89  &   96\\
         5 & 94 & 96.95  &   96\\ 
         6 & 94 & 97.16  &   98\\ 
         \hline
    \end{tabular}
    \caption{These values represent the similarity of solutions given from \eqref{withoutInteractionModel} and \eqref{withInteractionModel} over budgets varying from 5 to 95 with stepsizes of 5. For the 100 parcels, the number of parcels with the same protection status is recorded. The minimum, average, and median of this is computed and displayed. Case 1: $\{S_0, S_1\}$, Case 2: $\{S_2, S_3\}$, Case 3: $ \{S_4, S_5\}$, Case 4: $\{S_6, S_7\}$, Case 5: $\{S_0, S_1, S_2, S_3, S_4\}$, and Case 6: $\{S_5, S_6, S_7, S_0, S_1\}$ were all conducted.  }
    \label{tab:similar}
\end{table}

\end{document}